\documentclass[10pt]{article}
\usepackage{times,amsmath,amssymb,theorem,wrapfig,multicol,graphicx,color,hyperref}
\makeatletter
\def\section{\@startsection{section}{1}{0pt}{-3.25ex plus -1ex minus 
-.2ex}{1.5ex plus .2ex minus .3ex}{\normalfont\large\bf}}
\renewenvironment{abstract}
{}
{}
\newcommand{\parag}[1]{\vspace{-1ex}\paragraph{#1.}\!\!}
\newcommand{\defn}[1]{{\textit{\textbf{#1}}}}
\newcommand{\eg}{\emph{e.g.}}
\theoremheaderfont{\scshape}
\theorembodyfont{\normalfont\slshape}
\theoremstyle{plain}
\newtheorem{definition}{Definition}
\newtheorem{lemma}{Lemma}
\newtheorem{theorem}{Theorem}
\newenvironment{proof}{\begin{trivlist}\item{}\normalfont\vspace{-.8ex}\textit{Proof.}}{\hfill$\square$\vspace{-.5ex}\end{trivlist}}
\newcommand{\graphof}[1]{G(#1)}
\newcommand{\vars}{\mathcal{V}}
\newcommand{\atoms}{\mathcal{A}}
\newcommand{\peircepq}{((p\mkern-2mu\Rightarrow\mkern-2mu q)\mkern-2mu\Rightarrow\mkern-2mu p)\mkern-2mu\Rightarrow\mkern-2mu p}
\newcommand{\lift}[1]{\widehat{\mkern.2mu#1\mkern-.2mu}}
\newcommand{\restr}[2]{#1_{\restriction #2}}
\newcommand{\domrestr}[2]{#1^{\restriction #2}}
\newcommand{\putlabel}[3]{\put(#1,#2){\makebox(0,0){#3}}}
\newcommand{\putmathlabel}[3]{\put(#1,#2){\makebox(0,0){$#3$}}}
\newcommand{\putovercenteredlabel}[3]{\put(#1,#2){\makebox(0,0)[bc]{#3}}}
\newcommand{\circl}{\raisebox{-.69pt}[5.7pt][0pt]{\Large$\circ$}}
\newcommand{\bigbulle}{\raisebox{-.6pt}[5.8pt][0pt]{\Large$\bullet$}}
\newcommand{\bigcircl}{\raisebox{-1.3pt}[5.7pt][0pt]{\LARGE$\circ$}}
\newcommand{\bigsquarsymb}{\raisebox{0pt}[4.9pt][0pt]{\setlength\fboxrule{.53pt}\setlength\fboxsep{0pt}\framebox{\rule{1.8mm}{0mm}\rule{0mm}{1.8mm}}}}
\definecolor{grey}{gray}{.6}
\newcommand{\bigblacksquar}{\raisebox{.05ex}{\rule{1.8mm}{1.8mm}}}
\newcommand{\biggreysolidsquar}{\color{grey}\bigblacksquar}
\newcommand{\biggreysquarsymb}{\putlabel{0}{0}{\biggreysolidsquar}\putlabel{0}{0}{\bigsquarsymb}}
\newcommand{\biggreysquar}{\raisebox{-.05ex}{$\biggreysquarsymb$}}
\newcommand{\putbigcircl}[2]{\putlabel{#1}{#2}{\bigcircl}}
\newcommand{\putbigsquar}[2]{\putlabel{#1}{#2}{\biggreysquar}}
\newcommand{\putbigbulle}[2]{\putlabel{#1}{#2}{\bigbulle}}
\newcommand{\putbiglabelledbullet}[5]{\putbigbulle{#1}{#2}\put(#1,#2){\put(#4,#5){\makebox(0,0){#3}}}}
\newcommand{\putlabelledbullet}[5]{\put(#1,#2){\makebox(0,-.7){\large$\bullet$}\put(#4,#5){\makebox(0,0){#3}}}}
\newcommand{\putsmalllabelledbullet}[5]{\put(#1,#2){\makebox(0,-.7){$\bullet$}\put(#4,#5){\makebox(0,0){#3}}}}
\newcommand{\hdownleftmapsto}{\thicklines\put(0,0){\put(0,0){\vector(-1,-4){12.5}}}}
\newcommand{\hdownrightmapsto}{\thicklines\put(0,0){\put(0,0){\vector(1,-4){12.5}}}}
\newcommand{\thmgap}{\vspace*{-.5ex}}
\newcommand{\postthmgap}{\vspace*{-.2ex}}

\title{\vspace*{-2.8ex}\Large\bf Proofs Without Syntax}\author{\\[-2.8ex]\normalsize\sc Dominic J.~D.~Hughes \\[-.1ex]
\normalsize Stanford University}\date{}
\begin{document}
\maketitle
\vspace*{-2.8ex}

\begin{abstract}\noindent\normalfont\normalsize
\vspace*{-2ex}
\begin{quotation}\begin{center}\textsl{``[M]athematicians care no more for logic than logicians for mathematics.''}\hspace*{4ex}

\mbox{}\hspace*{48ex}Augustus de Morgan, 1868\end{center}\end{quotation}
\vspace*{-1.5ex}
\newlength{\indentlen}\setlength{\indentlen}{\parindent}
\begin{center}\begin{minipage}{5in}\hspace*{\indentlen}Proofs
are traditionally syntactic, inductively generated objects.  This
paper presents an abstract mathematical formulation of propositional
calculus (propositional logic) in which proofs are combinatorial
(graph-theoretic), rather than syntactic.  It defines a
\emph{combinatorial proof} of a proposition $\phi$ as a graph
homomorphism $h:C\to \graphof{\phi}$, where $\graphof{\phi}$ is a
graph associated with $\phi\mkern2mu$ and $C$ is a coloured graph.  The main
theorem is soundness and completeness: $\,\phi$ is true iff there
exists a combinatorial proof $h:C\to \graphof{\phi}$.\vspace*{2ex}
\end{minipage}\end{center}\end{abstract}

\begin{multicols}{2}
\thispagestyle{empty}\section{Introduction}\vspace*{1.5ex}
In 1868, de Morgan lamented the rift between \mbox{mathematics} and logic
\cite{deM68}: ``\textsl{[M]athematicians care no more for logic than
logicians for mathematics}.''  The dry syntactic manipulations of
formal logic can be off-putting to mathematicians accustomed to
beautiful symmetries, geometries, and rich layers of structure.
Figure~1 shows a syntactic proof in a standard Hilbert system taught
to mathematics undergraduates \cite{Hil28,Joh87}.%
\newcommand{\postfiggap}{-1.5ex}\newcommand{\axiomsgap}{-1.7ex}\newcommand{\preproofgap}{-1.7ex}\newcommand{\lineruleheight}{.83em}%
\newlength{\tw}\setlength{\tw}{\textwidth}\addtolength{\tw}{-\arrayrulewidth}\addtolength{\tw}{-2\tabcolsep}%
\newcommand{\twocolboxedfigt}[3]{\begin{figure*}[t]\vspace{-1.5ex}\begin{flushleft}\rule{5mm}{0.4pt}\hspace*{0.3mm}\raisebox{-3mm}{\rule{0mm}{6mm}}
\raisebox{-0.8mm}{\bf \,#2\,}{\rule{0mm}{6mm}}
\hrulefill{\rule{0mm}{6mm}}
\raisebox{-0.8mm}{\bf \,$\mkern1mu$\emph{\normalsize Figure #1}$\mkern-6mu$}%
\raisebox{-3mm}{\rule{0mm}{6mm}}\hspace*{1.3mm}\rule{5mm}{0.4pt}\raisebox{-3mm}{\rule{0mm}{6mm}}%
\vspace*{-3.4mm}
\begin{tabular}{@{}|c|@{}}\begin{minipage}[b]{\tw}\vspace*{5mm}\begingroup
\begin{multicols}{2}#3\end{multicols}\protect\rule{0mm}{1mm}\endgroup\end{minipage}\\ 
\hline\end{tabular}\end{flushleft}\vspace*{\postfiggap}\end{figure*}}
\twocolboxedfigt{1\label{hilbertproof}}
{A syntactic proof of $\,\peircepq\,$ in a standard Hilbert system}
{\newcommand\nn[1]{(#1\mkern-4mu\Rightarrow\mkern-7mu\bot\mkern-2mu)}%
\newcommand\fbot{\bot\mkern-1mu}%
\newcommand{\rr}{{\mkern1.5mu\Rightarrow\mkern0.5mu}}%
\newcommand\mm[2]{(m^{#1}_{#2})\mkern-4mu}%
{\small\vspace*{-3ex}Below is a proof of Peirce's law
$\,\peircepq\;$ in a standard Hilbert formulation of propositional
logic, taught to mathematics undergraduates \cite{Joh87}, with axiom
schemata\vspace*{\axiomsgap}\begin{center}\begin{math}
\begin{array}{rl}
(a) & x\Rightarrow (y\Rightarrow x) \\
(b) & (x\Rightarrow(y\Rightarrow z))\Rightarrow((x\Rightarrow y)\Rightarrow(x\Rightarrow z)) \\
(c) & ((x\Rightarrow \bot)\Rightarrow \bot)\Rightarrow x\\[\axiomsgap]
\end{array}
\end{math}\end{center}
and where $(m^i_j)$ marks \emph{modus ponens} with hypotheses numbered
$i$ and $j$.  \rule{0mm}{.7em}Hilbert systems tend to emphasise the elegance of the
\rule{0mm}{\lineruleheight}schemata (just $(a)$--$(c)$ suffice) over the elegance of the
proofs \rule{0mm}{\lineruleheight}generated by the schemata.
(Note: there may exist a shorter proof of Peirce's law in this\rule{0mm}{\lineruleheight}
system.)}
\vspace*{\preproofgap}{\tiny{\begin{displaymath}
\hspace{-5ex}\begin{array}{@{}r@{\;\;}l@{\;\;\;}l}
 1\mkern-2mu & (c)         & \nn{\nn{q}}\rr q \\[.5ex]
 2\mkern-2mu & (a)         & (\nn{\nn{q}}\rr q)\rr (\fbot\rr (\nn{\nn{q}}\rr q)) \\[.5ex]
 3\mkern-2mu & \mm{ 1}{ 2} & \fbot\rr (\nn{\nn{q}}\rr q) \\[.5ex]
 4\mkern-2mu & (b)         & (\fbot\rr (\nn{\nn{q}}\rr q))\rr ((\fbot\rr \nn{\nn{q}})\rr (\fbot\rr q)) \\[.5ex]
 5\mkern-2mu & \mm{ 3}{ 4} & (\fbot\rr \nn{\nn{q}})\rr (\fbot\rr q) \\[.5ex]
 6\mkern-2mu & (a)         & \fbot\rr \nn{\nn{q}} \\[.5ex]
 7\mkern-2mu & \mm{ 6}{ 5} & \fbot\rr q \\[.5ex]
 8\mkern-2mu & (a)         & (\fbot\rr q)\rr (p\rr (\fbot\rr q)) \\[.5ex]
 9\mkern-2mu & \mm{ 7}{ 8} & p\rr (\fbot\rr q) \\[.5ex]
10\mkern-2mu & (b)         & (p\rr (\fbot\rr q))\rr (\nn{p}\rr (p\rr q)) \\[.5ex]
11\mkern-2mu & \mm{ 9}{10} & \nn{p}\rr (p\rr q) \\[.5ex]
12\mkern-2mu & (a)         & ((p\rr q)\rr p)\rr (\nn{p}\rr ((p\rr q)\rr p)) \\[.5ex]
13\mkern-2mu & (b)         & (\nn{p}\rr ((p\rr q)\rr p))\rr ((\nn{p}\rr (p\rr q))\rr (\nn{p}\rr p)) \\[.5ex]
14\mkern-2mu & (a)         & ((\nn{p}\rr ((p\rr q)\rr p))\rr ((\nn{p}\rr (p\rr q))\rr (\nn{p}\rr p)))\Rightarrow\\&& (((p\rr q)\rr p)\rr ((\nn{p}\rr ((p\rr q)\rr p))\Rightarrow\\&& ((\nn{p}\rr (p\rr q))\rr (\nn{p}\rr p)))) \\[.5ex]
15\mkern-2mu & \mm{13}{14} & ((p\rr q)\rr p)\rr ((\nn{p}\rr ((p\rr q)\rr p))\Rightarrow\\&& ((\nn{p}\rr (p\rr q))\rr (\nn{p}\rr p))) \\[.5ex]
16\mkern-2mu & (b)         & (((p\rr q)\rr p)\rr ((\nn{p}\rr ((p\rr q)\rr p))\rr ((\nn{p}\rr (p\rr q))\Rightarrow\\&& (\nn{p}\rr p))))\rr ((((p\rr q)\rr p)\rr (\nn{p}\rr ((p\rr q)\rr p)))\Rightarrow\\&& (((p\rr q)\rr p)\rr ((\nn{p}\rr (p\rr q))\rr (\nn{p}\rr p)))) \\[.5ex]
17\mkern-2mu & \mm{15}{16} & (((p\rr q)\rr p)\rr (\nn{p}\rr ((p\rr q)\rr p)))\Rightarrow\\&& (((p\rr q)\rr p)\rr ((\nn{p}\rr (p\rr q))\rr (\nn{p}\rr p))) \\[5ex]
\end{array}\hspace{-5ex}\vspace*{-14ex}
\end{displaymath}
\begin{displaymath}
\hspace{-5ex}\begin{array}{@{}r@{\;\;}l@{\;\;\;}l}
\\[-2.5ex]18\mkern-2mu & \mm{12}{17} & ((p\rr q)\rr p)\rr ((\nn{p}\rr (p\rr q))\rr (\nn{p}\rr p)) \\[.65ex]
19\mkern-2mu & (b)         & (((p\rr q)\rr p)\rr ((\nn{p}\rr (p\rr q))\rr (\nn{p}\rr p)))\Rightarrow\\&& ((((p\rr q)\rr p)\rr (\nn{p}\rr (p\rr q)))\rr (((p\rr q)\rr p)\rr (\nn{p}\rr p)))\hspace*{-10ex}\\[.65ex]
20\mkern-2mu & \mm{18}{19} & (((p\rr q)\rr p)\rr (\nn{p}\rr (p\rr q)))\rr (((p\rr q)\rr p)\rr (\nn{p}\rr p)) \\[.65ex]
21\mkern-2mu & (a)         & (\nn{p}\rr (p\rr q))\rr (((p\rr q)\rr p)\rr (\nn{p}\rr (p\rr q))) \\[.65ex]
22\mkern-2mu & (a)         & ((((p\rr q)\rr p)\rr (\nn{p}\rr (p\rr q)))\rr (((p\rr q)\rr p)\rr (\nn{p}\rr p)))\hspace*{-10ex}\\&&\rr ((\nn{p}\rr (p\rr q))\rr ((((p\rr q)\rr p)\rr (\nn{p}\rr (p\rr q)))\Rightarrow\\&& (((p\rr q)\rr p)\rr (\nn{p}\rr p)))) \\[.65ex]
23\mkern-2mu & \mm{20}{22} & (\nn{p}\rr (p\rr q))\rr ((((p\rr q)\rr p)\rr (\nn{p}\rr (p\rr q)))\Rightarrow\\&& (((p\rr q)\rr p)\rr (\nn{p}\rr p))) \\[.65ex]
24\mkern-2mu & (b)         & ((\nn{p}\rr (p\rr q))\rr ((((p\rr q)\rr p)\rr (\nn{p}\rr (p\rr q)))\Rightarrow\\&& (((p\rr q)\rr p)\rr (\nn{p}\rr p))))\rr (((\nn{p}\rr (p\rr q))\Rightarrow\\&& (((p\rr q)\rr p)\rr (\nn{p}\rr (p\rr q))))\rr((\nn{p}\rr (p\rr q))\Rightarrow\\&& (((p\rr q)\rr p)\rr (\nn{p}\rr p)))) \\[.65ex]
25\mkern-2mu & \mm{23}{24} & ((\nn{p}\rr (p\rr q))\rr (((p\rr q)\rr p)\rr (\nn{p}\rr (p\rr q))))\Rightarrow\\&& ((\nn{p}\rr (p\rr q))\rr (((p\rr q)\rr p)\rr (\nn{p}\rr p))) \\[.65ex]
26\mkern-2mu & \mm{21}{25} & (\nn{p}\rr (p\rr q))\rr (((p\rr q)\rr p)\rr (\nn{p}\rr p)) \\[.65ex]
27\mkern-2mu & \mm{11}{26} & ((p\rr q)\rr p)\rr (\nn{p}\rr p) \\[.65ex]
28\mkern-2mu & (a)         & \nn{p}\rr ((\nn{p}\rr \nn{p})\rr \nn{p}) \\[.65ex]
29\mkern-2mu & (b)         & (\nn{p}\rr ((\nn{p}\rr \nn{p})\rr \nn{p}))\Rightarrow\\&& ((\nn{p}\rr (\nn{p}\rr \nn{p}))\rr (\nn{p}\rr \nn{p})) \\[.65ex]
30\mkern-2mu & \mm{28}{29} & (\nn{p}\rr (\nn{p}\rr \nn{p}))\rr (\nn{p}\rr \nn{p}) \\[.65ex]
31\mkern-2mu & (a)         & \nn{p}\rr (\nn{p}\rr \nn{p}) \\[.65ex]
32\mkern-2mu & \mm{31}{30} & \nn{p}\rr \nn{p} \\[.65ex]
33\mkern-2mu & (b)         & (\nn{p}\rr \nn{p})\rr ((\nn{p}\rr p)\rr \nn{\nn{p}}) \\[.65ex]
34\mkern-2mu & \mm{32}{33} & (\nn{p}\rr p)\rr \nn{\nn{p}} \\[.65ex]
35\mkern-2mu & (c)         & \nn{\nn{p}}\rr p \\[.65ex]
36\mkern-2mu & (a)         & (\nn{\nn{p}}\rr p)\rr ((\nn{p}\rr p)\rr (\nn{\nn{p}}\rr p)) \\[.65ex]
37\mkern-2mu & \mm{35}{36} & (\nn{p}\rr p)\rr (\nn{\nn{p}}\rr p) \\[.65ex]
38\mkern-2mu & (b)         & ((\nn{p}\rr p)\rr (\nn{\nn{p}}\rr p))\Rightarrow\\&& (((\nn{p}\rr p)\rr \nn{\nn{p}})\rr ((\nn{p}\rr p)\rr p)) \\[.65ex]
39\mkern-2mu & \mm{37}{38} & ((\nn{p}\rr p)\rr \nn{\nn{p}})\rr ((\nn{p}\rr p)\rr p) \\[.65ex]
40\mkern-2mu & \mm{34}{39} & (\nn{p}\rr p)\rr p \\[.65ex]
41\mkern-2mu & (a)         & ((\nn{p}\rr p)\rr p)\rr (((p\rr q)\rr p)\rr ((\nn{p}\rr p)\rr p)) \\[.65ex]
42\mkern-2mu & \mm{40}{41} & ((p\rr q)\rr p)\rr ((\nn{p}\rr p)\rr p) \\[.65ex]
43\mkern-2mu & (b)         & (((p\rr q)\rr p)\rr ((\nn{p}\rr p)\rr p))\Rightarrow\\&& ((((p\rr q)\rr p)\rr (\nn{p}\rr p))\rr (((p\rr q)\rr p)\rr p)) \\[.65ex]
44\mkern-2mu & \mm{42}{43} & (((p\rr q)\rr p)\rr (\nn{p}\rr p))\rr (((p\rr q)\rr p)\rr p) \\[.65ex]
45\mkern-2mu & \mm{27}{44} & ((p\rr q)\rr p)\rr p\\[-16ex]
\end{array}\hspace{-5ex}\end{displaymath}}}\mbox{\strut}} Although
the system itself is elegant (just three axiom schemata suffice),
the syntactic proofs generated in it need not be.  Other syntactic
systems include \cite{Fr1879,Gen35}.

This paper presents an abstract mathematical formulation of
propositional calculus (propositional logic) in which proofs are
combinatorial (graph-theoretic), rather than syntactic.  It defines a
\emph{combinatorial proof\,} of a proposition $\phi\,$ as a
graph \mbox{homomorphism} \mbox{$h:C\to \graphof{\phi}$}, where $\graphof{\phi}$ is a graph
associated with $\phi$ and $C$ is a coloured graph.  
For example, if $\;\phi\:=\:\peircepq\;$ 
then $\graphof{\phi}$ is:
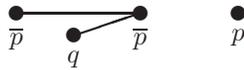
\begin{center}\begin{picture}(0,19)(0,-14)\setlength{\unitlength}{.7pt}
\put(-60,0){\thicklines
\putbiglabelledbullet{0}{0}{$\overline{p}$}{0}{-14}
\put(0,0){\line(1,0){67}}
\putbiglabelledbullet{31}{-12}{$q$}{0}{-13}
\putbiglabelledbullet{67}{0}{$\overline{p}$}{0}{-14}
\put(31,-12){\line(3,1){35}}
\put(120,0){\putbiglabelledbullet{0}{0}{$p$}{0}{-13}}}\end{picture}\end{center}
A combinatorial proof $h:C\to \graphof{\phi}$ of $\phi$ is shown below:
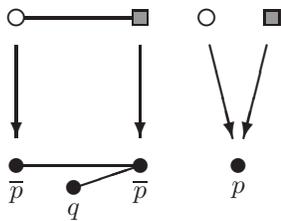
\begin{center}\begin{picture}(0,75)(0,-16)\put(0,0){\setlength{\unitlength}{.7pt}
\put(-60,0){\thicklines
\putbiglabelledbullet{0}{0}{$\overline{p}$}{0}{-14}
\put(0,0){\line(1,0){67}}
\putbiglabelledbullet{31}{-12}{$q$}{0}{-13}
\putbiglabelledbullet{67}{0}{$\overline{p}$}{0}{-14}
\put(31,-12){\line(3,1){35}}
\put(120,0){\putbiglabelledbullet{0}{0}{$p$}{0}{-13}}}}\setlength{\unitlength}{.7pt}
\put(-60,0){\thicklines
\putbigcircl{0}{80}
\put(0,0){\line(1,0){68}}
\put(4.5,80){\line(1,0){58.5}}
\put(0,65){\thicklines\vector(0,-1){50}}
\putbigsquar{67}{80}
\put(67,65){\thicklines\vector(0,-1){50}}
\put(120,0){\putbigcircl{-17}{80}
\putbigsquar{18}{80}
\put(-15.5,65){\hdownrightmapsto}
\put(15.5,65){\hdownleftmapsto}}}\end{picture}\end{center}\label{peirceeg}
The upper graph $C$ has two colours (white $\bigcircl$ and grey
\raisebox{.7ex}{$\mkern9mu\biggreysquar\mkern9mu$}), and the arrows
define $h$.  The same proposition is proved syntactically in Figure~1.

The main theorem of the paper is soundness and completeness:
\begin{center}\vspace{-.2ex}\textsl{A proposition is true iff it has a combinatorial proof.}\begin{picture}(0,0)
\put(29,190){\rotatebox[origin=l]{270}{\LARGE\color[rgb]{.35,.35,.35}{To appear in \textsl{Annals of Mathematics}.  Submitted 20 Aug 2004, accepted 9 Sep 2005.}}}\end{picture}\vspace{-.2ex}\end{center}
As with conventional syntactic soundness and completeness, this
theorem match\-es a universal quantification with an existential one:
a proposition $\phi$ is true if it evaluates to $1\,$ \emph{for all}
$\,0/1$ assignments of its variables, and $\phi$ is provable if
\emph{there exists} a proof of $\phi$.  However, where conventional
completeness provides an inductively generated \emph{syntactic}
witness (\eg\ Figure~1), this theorem provides an abstract
\emph{mathematical} witness for every true proposition (\eg\ the
homomorphism $h$ drawn above).

Just three conditions suffice for soundness and completeness: a graph
homomorphism $h:C\to G(\phi)$ is a combinatorial proof of $\phi$ if
(1) $C$ is a suitable coloured graph, (2) the image of each colour
class is labelled appropriately, and (3) $h$ is a \emph{skew
fibration}, a lax form of graph fibration.  Each condition can be
checked in polynomial time, so combinatorial proofs constitute a
formal \emph{proof system} \cite{CR79}.

\vspace*{-.2ex}\parag{Acknowledgements}\hspace*{-.5ex}
Nil Demir\c{c}ubuk, Vaughan Pratt, Ju\-lien Basch, Ra\-jat
Bhat\-ta\-char\-jee, Rob van Glab\-beek, Sol Feferman, Grisha Mints
and Don Knuth.  Stanford grant {\small 1DMA644}.

\vspace{-.1ex}\section{Notation and terminology}\label{notation}

\vspace*{2.4ex}
\parag{Graphs}
An \defn{edge} on a set $V$ is a two-element subset of $V\!\!$.\, A
\defn{graph} $(V,E)$ is a finite set $V\!\!$ of \defn{vertices}
and a set $E$ of edges on $V\!\!$.  Write $V(G)$ and $E(G)$ for the
vertex set and edge set of a graph $G$, respectively, and $vw$ for
$\{v,w\}$.  The \defn{complement} of $(V,E)$ is the graph
$(V,E^{\mkern1mu\mathsf{c}})$ with $vw\mkern-3mu\in\mkern-3mu
E^{\mkern1mu\mathsf{c}}$ iff $vw\mkern-3mu\not\in\mkern-4mu E$.  A
graph $(V,E)$ is \defn{coloured}\/ if $V$ carries an equivalence
relation $\sim$ such that $v\mkern-3mu\sim\mkern-3mu w$ only if
$vw\mkern-2mu\not\in\mkern-2mu E$; each equivalence class is a
\defn{colour class}.  Given a set $L$, a graph is
\defn{\mbox{$L$-labelled}} if every vertex has an \mbox{element} of
$L$ associated with it, its \defn{label}.  Let $G\mkern-3mu
=\mkern-3mu (V\mkern-4mu,\mkern-2mu E)$ and $G'\mkern-5mu =\mkern-3mu
(V'\mkern-5mu,\mkern-2mu E')$ be graphs.  A \defn{homomorphism}
$\mkern2mu h:G\to G'\mkern2mu$ is a function $\mkern2mu h:\mkern-2mu
V\to V'$ such that $vw\in E$ implies $\mkern3mu h(v)h(w)\in E'$.  If
$V$ and $V'$ are disjoint, the \defn{union}
$G\mkern-.1mu\vee\mkern-.1mu G'$ is $(V\mkern-.2mu\cup\mkern-.05mu
V'\mkern-5mu,\mkern1mu E\mkern1.5mu\cup\mkern-.3mu E')$ and the
\defn{join} $\mkern2muG\mkern-1mu\wedge G'$ is $(V\cup
V'\!\!,\mkern2mu E\cup E'\cup\{\mkern2mu vv'\mkern-2mu:v\!\in\!
V\!,\:v'\!\in\! V'\})$; colourings or labellings of $G$ and $G'$ are
inherited.  A graph $(V,E)$ is a \defn{cograph} \cite{CLS81} if $V$ is
non-empty and for any distinct $v,w,x,y\!\in\!  V\mkern-1mu$, the
restriction of $E$ to edges on $\{v,w,x,y\}$ is not $\{vw,wx,xy\}$.  A
set $W\mkern-2mu\subseteq\mkern-2mu V\mkern-1mu$ \defn{induces a
matching} if it is non-empty and for all $w\mkern-2mu\in\mkern-3mu W$
there is a unique $w'\mkern-3mu\in\mkern-3mu W$ such that
$ww'\mkern-2mu\in\mkern-1mu E$.

\vspace*{-5.5ex}\parag{Propositions}\hspace*{-1.1ex}
Fix a set $\vars\mkern-1mu$ of \defn{variables}.
\hspace*{-.1ex}A \defn{proposition} is any expression generated freely from variables
by the binary operations \defn{and} $\wedge$,
\defn{or} $\vee$, and \defn{implies} $\Rightarrow$,
the unary operation
\defn{not} $\neg$,
and the constants (nullary operations) \defn{true} $1$ and
\defn{false} $0$.
A \defn{valuation}\label{valuation} is a function $f:\vars\to\{0,1\}$.
Write $\hat{f}$ for the extension of a valuation $f$ to propositions
defined by $\hat{f}(0)\mkern-2mu=\mkern-2mu 0$,
$\hat{f}(1)\mkern-2mu=\mkern-2mu 1$,
$\hat{f}(\neg\phi)\mkern-2mu=\mkern-2mu1\mkern-4mu-\mkern-4mu\hat{f}(\phi)$,
$\hat{f}(\phi\mkern-3mu\wedge\mkern-3mu
\rho)\mkern-3mu=\mkern-3mu\min\{\hat{f}(\phi),\mkern-2mu\hat{f}(\rho)\}$,
$\hat{f}(\phi\mkern-1mu\vee\mkern-3mu\rho)\mkern-3mu=\mkern-3mu\max\{\hat{f}(\phi),\mkern-2mu\hat{f}(\rho)\}$,
$\hat{f}(\phi\mkern-1mu\Rightarrow\!\rho)\mkern-3mu=\mkern-3mu
\hat{f}((\neg\phi)\mkern-2mu\vee\mkern-3mu\rho)$.
A proposition $\phi$ is \defn{true} if
$\hat{f}(\phi)\mkern-3mu=\mkern-3mu 1$ for all valuations $f$.
Variables $p\in\vars$ and their negations
$\overline{p}\mkern-2mu=\mkern-3mu\neg p$ are \defn{literals}; \;$p$
and $\overline{p}$ are \defn{dual}, as are $0$ and $1$.  An
\defn{atom} is a literal or constant, and $\atoms$ denotes the set of
atoms.%
\section{Combinatorial proofs}\label{comb-proofs}

\mbox{}\vspace*{-\baselineskip}\vspace*{-\parskip}

\newcommand{\width}{20}\newcommand{\halfwidth}{10}\begin{wrapfigure}{r}{.63in}
\begin{center}\begin{picture}(0,0)(-1,-9)%
\thicklines\put(-\halfwidth,-\halfwidth){\putbiglabelledbullet{0}{\width}{$p$}{-9}{0}%
\putbiglabelledbullet{0}{0}{$\overline{q}$}{-9}{0}\putbiglabelledbullet{\width}{\width}{$0$}{8}{0}%
\putbiglabelledbullet{\width}{0}{$p$}{8}{0}\put(0,\width){\line(1,0){\width}}%
\put(0,\width){\line(1,-1){\width}}\put(0,0){\line(1,0){\width}}%
\put(0,0){\line(1,1){\width}}}
\end{picture}\end{center}\vspace*{-2ex}\end{wrapfigure}
\noindent Given an $\atoms$-labelled graph $G$, define $\neg G$ as the
result of complementing $G$ and every label of $G$.  For example, if
$G$ is the graph shown right, then $\neg G$ is the graph below left.
Define $\:G\mkern-3mu\Rightarrow\mkern-3mu G'\mkern2mu=\mkern2mu 
(\neg G)\mkern-2mu\vee\mkern-2mu G'.$
Identify each atom $a$ with a single vertex labelled $a$; thus,
having defined operations
$\neg$, $\vee$,
\mbox{\hspace{20ex}\hfill}\vspace*{-\baselineskip}\vspace*{-\parskip}
 
\begin{wrapfigure}{l}{.6in}
\begin{center}\begin{picture}(0,0)(11,0)
\thicklines\putbiglabelledbullet{0}{\halfwidth}{$\overline{p}$}{-8}{0}
\putbiglabelledbullet{0}{-\halfwidth}{$q$}{-8}{0}
\putbiglabelledbullet{\width}{\halfwidth}{$1$}{8}{0}
\putbiglabelledbullet{\width}{-\halfwidth}{$\overline{p}$}{8}{0}
\put(0,\halfwidth){\line(0,-1){\width}}
\put(\width,\halfwidth){\line(0,-1){\width}}
\end{picture}\end{center}
\end{wrapfigure}\noindent 
$\wedge$ and $\Rightarrow$ on $\atoms$-labelled graphs,
every proposition $\phi$
determines an $\atoms$-labelled graph, denoted $\graphof{\phi}$.
For example,
$G\big((p\vee\mkern-3mu \neg q)\wedge(0\vee p)\big)$
is above right,
$G\big((q\wedge \neg p)\vee(1\wedge \neg p)\big)$ 
is left,
and $G\big(\peircepq\big)$ is 
in the Introduction.

\begin{wrapfigure}{r}{.94in}
\setlength{\unitlength}{.8pt}
\begin{center}
\begin{picture}(0,19)(10,16)
\thicklines
\put(0,50){\putlabel{-15}{1}{$\forall v$}\putlabel{35}{7}{$\lift{w}$}\put(-4,0){\line(6,1){31}}\put(-15,-9)
{\vector(0,-1){30}}\put(35,-2){\vector(0,-1){30}}\putlabel{11}{10}{$\exists$}
}
\putlabel{-15}{0}{$h(v)$}
\putlabel{36}{8}{$h(\lift{w})$}
\putlabel{21}{-10}{$w$}
\put(-3,-3){\line(5,-1){17}}
\put(-3,1){\line(5,1){23}}
\putlabel{4}{-11}{$\forall$}
\end{picture}
\end{center}
\end{wrapfigure}
A colouring is \defn{nice} if every colour class has at most two
vertices and no union of two-vertex colour classes induces a matching.
A graph homomorphism $h\mkern-1mu:\mkern-1mu G\mkern-2mu\to\mkern-2mu
G'$ is a \defn{skew fibration}\label{skewfib} (see figure
right\label{skew-fig}) if for all $v\mkern-2mu\in\mkern-2mu V(G)$ and
$h(v)w\mkern-2mu\in\mkern-2mu E(G')$ there exists
$v\lift{w}\mkern-2mu\in\mkern-2mu E(G)$ with
$h(\lift{w})w\mkern-2mu\not\in\mkern-2mu E(G')$.  Given a graph
homomorphism $h:G\to G'$ with $G'$ an $\atoms$-labelled graph, a
vertex $v\in V(G)$ is \defn{\mbox{axiomatic}} if $h(v)$ is labelled
$1$, and a pair $\{v,w\}\subseteq V(G)$ is axiomatic if $h(v)$ and
$h(w)$ are labelled by dual literals.
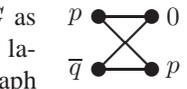
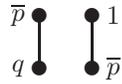
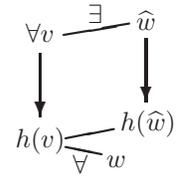
\begin{definition}\label{combproof}
A \defn{combinatorial proof} of a proposition $\phi$ is a skew
fibration \mbox{$h:C\to \graphof{\phi}$} from a nicely coloured
\mbox{cograph} $C$ to the graph $G(\phi)$ of $\phi$, such that every
colour class of $C$ is axiomatic.
\end{definition}%
\postthmgap\vspace*{-.5ex}
A combinatorial proof of 
$\,\peircepq\,$ is shown in the Introduction.
The reader may find it instructive to consider why $p\wedge\! \neg p\,$
has no combinatorial proof.
\thmgap\vspace*{-.5ex}
\begin{theorem}[Soundness and Completeness]\label{thm}\mbox{}\\
A proposition is true iff it has a combinatorial proof.
\end{theorem}
\postthmgap\vspace*{-.2ex}
Section~\ref{combprops} reformulates this theorem in terms of
combinatorial (non-syntactic, non-inductive) notions of
\emph{proposition} and \emph{truth}.  
Section~\ref{proof} proves the reformulated theorem.

\parag{Notes}\label{notes}
The map $\phi\mapsto G(\phi)$ is based on a well understood
translation of a boolean formula into a graph \cite{CLS81}, and (up to
standard graph isomorphism\footnote{Graphs $(V,E)$ and $(V',E')$ are
isomorphic if there exists a bijection $h:V\to V'$ with $vw\in E$ iff
$h(v)h(w)\in E'$.\label{iso-footnote}}) represents propositions modulo
associativity and commutativity of $\wedge$ and $\vee$, double
negation $\neg\neg\phi=\phi$, de Morgan duality
$\neg(\phi\wedge\rho)=(\neg\phi)\vee (\neg
\rho)$ and $\neg(\phi\vee\rho)=(\neg\phi)\wedge (\neg \rho)$, and
$\phi\Rightarrow \rho\:=\:(\neg \phi)\vee \rho$.  Perhaps the
earliest graphical representation of propositions is due to Peirce
\cite[vol.\,4:2]{Pei58}, dating from the late 1800s.

A skew fibration is a lax notion of graph fibration.  A graph homomorphism
$h:G\to G'$ is a \defn{graph fibration} (see \eg\
\cite{BV02}) if for all $v\in V(G)$ and
$h(v)w\in E(G')$ there is a unique $v\lift{w}\in E(G)$ with
$h(\lift{w})\!=\! w\,$.\footnote{This is simply a convenient
restatement of the familiar notions of fibration in topology
\cite{Whi78} and category theory \cite{Gro59,Gra66}: a graph
homomorphism is a graph fibration \emph{iff} it satisfies the homotopy
lifting property (when viewed as a continuous map by identifying each
edge with a copy of the unit interval) \emph{iff} it has all requisite
cartesian liftings (when viewed as a functor by identifying each graph
with its path category).} The definition of skew fibration drops
uniqueness and relaxes $h(\lift{w})\!=\!w$ to `skewness'
$h(\lift{w})w\not\in E(G')$.

Combinatorial proofs constitute a formal \emph{proof system}
\cite{CR79} since correctness can be checked in polynomial time.\footnote{The
skew fibration and axiomatic conditions are clearly
polynomial.  Checking that a graph $G$ is a cograph is
polynomial by constructing its modular decomposition tree $T(G)\,$
\cite{BLS99}, and checking that $G$ is nicely coloured is a simple
breadth-first search on $T(G)$.}  There is a polynomial-time
computable function taking a propositional sequent calculus proof of
$\phi$ with $n\ge 0$ cut rules \cite{Gen35} to a combinatorial proof
of $\phi$ \defn{with $n$ cuts}: a combinatorial proof of
$\phi\vee(\theta_1\wedge\neg\theta_1)\vee\cdots\vee(\theta_n\wedge\neg\theta_n)$
for propositions $\theta_i$.

In the example of a combinatorial proof drawn in the Introduction,
observe that the image of the colour class
$\:\circl\hspace*{1ex}\circl\:$ under $h$ is
\begin{picture}(30,0)(0,-7)\putlabelledbullet{8}{0}{\small$\overline{p}$}{0}{-7.3}
\putlabelledbullet{22}{0}{\small$p$}{0}{-7.5}\end{picture}.
Think of the colour class as actively pairing an occurrence of a
variable $p$ with an occurrence of its dual $\overline{p}$.  The idea
of pairing dual variable occurrences has arisen in the study of
various forms of syntax, such as closed categories
\cite{KM71}, contraction-free predicate calculus
\cite{KW84} and linear logic \cite{Gir87}.
Combinatorial proofs relate only superficially to the
connection/matrix method \cite{Dav71,Bib74,And81}; the latter fails to
provide a proof system \cite{CR79}.

A partially combinatorial notion of proof for classical logic, called
a \emph{proof net}, was presented in \cite{Gir91}, though promptly
dismissed by the author as overly syntactic: a proof net of a
proposition $\phi$ has an underlying syntax tree containing not only
$\wedge$'s and $\vee$'s from $\phi$, but also auxiliary syntactic
connectives which are not even boolean operations (\emph{contraction}
and \emph{weakening}).

Nicely coloured cographs with two vertices in every colour class
correspond to \emph{unlabelled chorded R\&B-cographs} \cite{Ret03}.
When labelled, the latter represent proof nets of mixed multiplicative
linear logic \cite{Gir87}.

\section{Combinatorial propositions and truth}\label{combprops}
\vspace{2ex}

A set
$W\!\mkern-3mu\subseteq\mkern-1mu\! V\mkern-2mu(G)\mkern-1mu$ is \defn{stable} if 
$vw\mkern-1mu\!\not\in\!\mkern-2mu E(G)\mkern-2mu$ for all $v\mkern-1mu,\mkern-2mu w\!\in\!\mkern-2mu W\!$.
A \defn{clause}
is a maximal stable set.
A clause of an $\atoms$-labelled graph is \defn{true} if it con\-tains
a $1$-labelled vertex or two vertices labelled by dual literals;
an $\atoms$-labelled graph is \defn{true} if its clauses are true.
For example,
{\setlength{\unitlength}{.85pt}\begin{picture}(36,12)(4,-7)\thicklines\putsmalllabelledbullet{10}{0}{\scriptsize$\overline{p}$}{0}{-7}\putsmalllabelledbullet{22}{0}{\scriptsize$p$}{0}{-7.5}
\putsmalllabelledbullet{34}{0}
{\scriptsize$1$}{0}{-7}
\put(22,0){\line(1,0){12}}\end{picture}}
\,(\,$=$ $G(p\mkern-2mu\Rightarrow\mkern-4mu(p\mkern-1mu\wedge\mkern-2mu 1))$\,) 
is true,
with true clauses
{\setlength{\unitlength}{.85pt}\begin{picture}(24,12)(4,-7)\thicklines\putsmalllabelledbullet{10}{0}{\scriptsize$\overline{p}$}{0}{-7}\putsmalllabelledbullet{22}{0}{\scriptsize$p$}{0}{-7.5}\end{picture}}
and 
{\setlength{\unitlength}{.85pt}\begin{picture}(36,12)(4,-7)\thicklines\putsmalllabelledbullet{10}{0}{\scriptsize$\overline{p}$}{0}{-7}\thicklines\putsmalllabelledbullet{34}{0}{\scriptsize$1$}{0}{-7}\end{picture}}.\vspace*{-1ex}
\begin{lemma}\label{true}
A proposition $\phi$ is true iff its
graph $G(\phi)$ is true.\postthmgap
\end{lemma}\begin{proof}
Exhaustively apply distributivity
$\theta\vee(\psi_1\wedge\psi_2)\to(\theta\vee\psi_1)\wedge(\theta\vee\psi_2)$
to $\phi$ modulo associativity and commutativity of $\wedge$ and
$\vee$, yielding a conjunction $\phi'$ of syntactic clauses (disjunctions of
atoms).  The lemma is immediate for $\phi'$
since $G(\phi')$ is a join of clauses,
and $G\big(\theta\vee(\psi_1\wedge\psi_2)\big)$
is true iff
$G\big((\theta\vee\psi_1)\wedge(\theta\vee\psi_2)\big)$
is true since
for non-empty graphs $G_1\mkern-2mu$ and $G_2$, a clause of $G_1\mkern-3mu\vee\mkern1mu G_2$
(resp.\ $G_1\mkern-1mu\wedge\mkern1mu G_2$) is a clause of $G_1\mkern-2mu$ and (resp.\ or) a clause of $G_2$.
\end{proof}A \defn{combinatorial proposition} is an $\atoms$-labelled cograph.
Since a graph is a cograph iff it is derivable from individual
vertices by union, join and complement
\cite[\S11.3]{BLS99}, 
the graph $G(\phi)$ of any syntactic proposition $\phi$ is a
combinatorial proposition; conversely every combinatorial proposition
is (isomorphic$^\text{\ref{iso-footnote}}$
to) $G(\phi)$ for some $\phi$.
\begin{definition}\label{abscombpf}\postthmgap\thmgap%
A \defn{combinatorial proof} of a combinatorial proposition $P$ is a
skew fibration $h:C\to P$ from a nicely coloured cograph $C$
whose colour classes are axiomatic.
\end{definition}\postthmgap
Thus a combinatorial proof of a syntactic proposition $\phi$ (Def.~\ref{combproof}) is a
combinatorial proof of $G(\phi)$ (Def.~\ref{abscombpf}).
By Lemma~\ref{true}, the following is equivalent to Theorem~\ref{thm}.
\vspace{-1ex}\begin{theorem}\label{absthm}\textsc{(Combinatorial Soundness and Completeness)}
\;A combinatorial proposition is true iff it has a combinatorial proof.
\end{theorem}\vspace{-1.5ex}

\section{Proof of Theorem~2}\label{proof}

\mbox{\rule{10ex}{0ex}}

\mbox{\rule{10ex}{0ex}}

\mbox{\rule{10ex}{0ex}}

\mbox{\rule{10ex}{0ex}}

\vspace{-4\baselineskip}\vspace{-4\parskip}
\begin{wrapfigure}{r}{1.14in}
\setlength{\unitlength}{2.5pt}
\newcommand{\putt}[3]{\putovercenteredlabel{#1}{#2}{#3}}
\newcommand{\vv}[3]{\put(0,0){\vector(#1,#2){#3}}\thicklines\put(0,0){\line(#1,#2){#3}}\thinlines}
\begin{picture}(0,16)(-1.7,2.95)
\putt{0}{0}{T\ref{thm}}
\putmathlabel{5}{1}{\leftrightarrow}
\putt{10}{0}{T\ref{absthm}}
\putt{5}{2.7}{\small\ref{true}}
\putt{20}{0}{T\ref{abs-comp}}
\putmathlabel{15}{1}{\leftarrow}
\putt{30}{0}{\ref{fusion-nice}}
\putmathlabel{25}{1}{\leftarrow}
\putt{0}{10}{\ref{clause}}
\putmathlabel{5}{11}{\rightarrow}
\putt{10}{10}{\ref{imgtrue}}
\putmathlabel{15}{11}{\rightarrow}
\putt{20}{10}{T\ref{thm-soundness}}
\putmathlabel{20}{6.2}{\downarrow}
\putmathlabel{25}{11}{\leftarrow}
\putt{30}{10}{\ref{tensor-true}}
\putmathlabel{15.3}{6.3}{\swarrow}
\putmathlabel{25.3}{6.3}{\swarrow}
\putt{0}{20}{\ref{skewfibs}}
\putmathlabel{5}{21}{\rightarrow}
\putmathlabel{0}{16.5}{\downarrow}
\putt{10}{20}{\ref{long-lemma}}
\putmathlabel{14.7}{16.6}{\searrow}
\putt{20}{20}{\ref{shallow}}
\putmathlabel{20}{16.5}{\downarrow}
\putt{30}{20}{\ref{fusion}}
\putmathlabel{25.6}{16.6}{\swarrow}
\end{picture}
\end{wrapfigure}
\noindent The diagram right shows the dependency between the Lemmas
(1--\ref{tensor-true}) and Theorems (T1--T4) in this paper.

Given a graph homomorphism $h:G\to G'$, an edge
$v\lift{w}\mkern-2mu\in\!E(G)$ is a \defn{skew lifting of $h(v)w\in
E(G')$ at $v$} if $h(\lift{w})w\mkern-2mu\not\in\mkern-2mu E(G')$.
Thus $h$ is a skew fibration iff every edge $h(v)w\in E(G')$ has a
skew lifting at $v$.

A graph $G$ is a \defn{subgraph} of $G'$, denoted $G\subseteq
G'$, if $V(G)\!\subseteq\! V(G')$ and $E(G)\!\subseteq\!  E(G')$.  
The subgraph $G[W]$ of $G$ \defn{induced by $\mkern1mu W\mkern-4mu\subseteq \mkern-3muV\mkern-2mu(G)$} is
$(W,\{\,vw\!\in\! E(G):v,\mkern-2mu w\!\in\! W\,\})$.
Let $h\mkern-1mu:\mkern-1mu G\mkern-2mu\to\mkern-2mu H$ be a graph
homomorphism and let $G'$ and $H'$ be induced subgraphs of $G$ and $H\!$,
respectively.
Write $h(G')$ 
for the induced 
subgraph $H\mkern4mu[\mkern1mu h(V(G'))\mkern1mu]$
and $h^{-1}(H')$ for the induced 
subgraph $G\mkern1mu[\mkern1mu h^{-1}(V(H'))\mkern1mu]$.
\mbox{Define} the \defn{restriction}
$\restr{h}{H'}\,:\,h^{-1}(H')\to H'\mkern2mu$ by $\mkern2mu\restr{h}{H'}(v)=h(v)$.

\begin{lemma}\label{skewfibs}
Let $\diamond\in\{\wedge,\vee\}$.  If $h:G\to H_1\diamond H_2$ is a skew fibration 
then both restrictions $\restr{h}{H_i}$ are skew fibrations.
\end{lemma}
\begin{proof}
We prove that if $v\lift{w}$ is a skew lifting of
$\restr{h}{H_i}(v)w=h(v)w\in E(H_i)$ at $v$ with respect to $h$, then
$h(\lift{w})\in H_i\,$; hence $v\lift{w}$ is a well-defined skew
lifting with respect to $\restr{h}{H_i}$.  Suppose $h(\lift{w})\in
H_j$ and $j\neq i$.  If $\diamond\!=\!\vee$, since $h$ is a
homomorphism, $h(v)h(\lift{w})$ is an edge between $H_1$ and $H_2$ in
$H_1\vee H_2$, a contradiction;
\,if $\mkern1mu\diamond\mkern-2mu=\!\wedge\:$, since $H_1\wedge H_2$ has all edges between
$H_1$ and $H_2$, $h(\lift{w})w$ is an edge, contradicting
$v\lift{w}$ being a skew lifting with respect to $h$.
\end{proof}
\begin{lemma}\label{long-lemma}
Let $\;h:(G_1\wedge G_2)\vee (H_1\vee H_2)\to (K_1\wedge K_2)\vee L$
be a skew fibration with $h(G_i)\subseteq K_i$ and $h(H_i)\subseteq
L$.  Then $h_i\mkern-2mu:\mkern-1mu G_i\mkern-2mu\vee\mkern-2mu
H_i\mkern-1mu\to\mkern-2mu K_i\mkern-2mu\vee\mkern-3mu L$ defined by
$h_i(v)\mkern-2mu =\mkern-2mu h(v)$ is a skew fibration.
\end{lemma}
\begin{proof}
Since a graph union $X_1\vee X_2$ has no edges between $X_1$ and
$X_2$, (a) if $k:X_1\vee X_2\to Y$ is a skew fibration, so also is
$\domrestr{k}{X_i}:X_i\to Y$ defined by $\domrestr{k}{X_i}(x)=k(x)$,
and (b) if $k_i:Z_i\to X_i$ is a skew fibration for $i=1,2$, so also
is $k_1\mkern-1.5mu\vee\mkern.5mu k_2:Z_1\vee Z_2\to X_1\vee X_2$ defined by
$(k_1\mkern-1.5mu\vee\mkern.5mu k_2)(z)=k_i(z)$ iff $z\in V(Z_i)$.  Since
$h_i\mkern4mu=\mkern6mu\restr{h}{K_i}\mkern2mu\vee\mkern5mu\domrestr{(\restr{h}{L})}{\mkern-1muH_i}$,
it is a skew fibration by (a), (b) and Lemma~\ref{skewfibs}.
\end{proof}
\begin{lemma}\label{clause}
If $h:G\to K$ is a skew fibration into a cograph $K$, then every
clause of $K$ contains a clause of $h(G)$.
\end{lemma}
\begin{proof}
By induction on the number of vertices in $K$.  The base case with $K$
a single vertex is immediate.  Otherwise $K=K_1\diamond K_2$ for
$\diamond\in\{\wedge,\vee\}$ and cographs $K_i$.  Let
$G_i=h^{-1}(K_i)$ and $h_i=\restr{h}{K_i}:G_i\to K_i$, a skew
fibration by Lemma~\ref{skewfibs}.  Let $C$ be a clause of $K$.  If
$\diamond=\wedge$ then $C$ is a clause of $K_j$ for $j=1$ or $2$; by
induction $C$ contains a clause $C'$ of $h_j(G_j)$, also a clause of
$h_1(G_1)\wedge h_2(G_2)=h(G)$.  If $\diamond=\vee$ then $C=C_1\cup
C_2$ for clauses $C_i$ of $K_i$; by induction $C_i$ contains a
clause $C'_i$ of $h_i(G_i)$, so $C$ contains the clause $C'_1\cup
C'_2$ of $h_1(G_1)\vee h_2(G_2)=h(G)$.
\end{proof}
\begin{lemma}\label{imgtrue}
Let $h:G\to P$ be a skew fibration into a combinatorial proposition
$P$.  If $h(G)$ is true then $P$ is true.
\end{lemma}
\begin{proof}
Lemma~\ref{clause} and the definition of \emph{true}.	
\end{proof}
The \defn{empty} graph is the graph with no vertices.  A graph is
\defn{disconnected} if it is a union of non-empty graphs, and
\defn{connected} otherwise.  A \defn{component} is a maximal non-empty connected
subgraph.  A graph homomorphism $h:G\to H$ is \defn{shallow} if
$h^{-1}(K)$ has at most one component for every component $K$ of $H$.
\vspace{-1.5ex}
\begin{lemma}\label{shallow}
For any combinatorial proof $\mkern1mu h:G\to P$ there exists a shallow
combinatorial proof $\mkern1mu h':G\to P'$ such that $P$ is true iff $P'$
is true.
\end{lemma}\vspace{-2.5ex}
\begin{proof}
Let $G_1,\ldots, G_n$ be the components of $G$, and let $P'$ be the
union of $n$ copies of $P$ defined by $V(P')=V(P)\times
\{1,\ldots,n\}$ and $\langle v,i\rangle\langle w,j\rangle\in E(P')$
iff $vw\in E(P)$ and $i=j$, and the label of $\langle v,i\rangle$ in
$P'$ equal to the label of $v$ in $P$.  Define $h':G\to P'$ on $v\in
V(G_i)$ by $h'(v)=\langle h(v),i\rangle$.  Since $P'$ is a union of
copies of $P\mkern-2mu$, it is true iff $P$ is true (every clause of
$P'$ contains a clause of $P$; conversely the union of $n$ copies of a
clause of $P$ is a clause of $P'$), and $h'$ is a combinatorial proof
(skew liftings copied from $h$).
\end{proof}
A subgraph $G'$ of $G$ is a \defn{portion} of $G$ if $G=G'\vee G''$
for some $G''$.  A \defn{fusion} of graphs $G$ and $H$ is any graph
obtained from $G\vee H$ by selecting portions $G'$ of $G$ and $H'$ of
$H$ and adding edges between every vertex of $G'$ and every vertex of
$H'$.  Union and join are extremal cases of fusion: union with $G',H'$
empty; join with $G'\!=\!G$, $H'\!=\!H$.  On coloured graphs, fusion
does not reduce to union and join: the coloured cograph
${\setlength{\unitlength}{1.4pt}\begin{picture}(47,0)(4,-2)\thicklines\putbigcircl{10}{0}\putbigcircl{20}{0}
\putbigsquar{33}{0}\putbigsquar{44}{0}\put(22,0){\line(1,0){8.9}}\end{picture}}$
is a fusion of
${\setlength{\unitlength}{1.4pt}\begin{picture}(20,0)(5,-2)\thicklines\putbigcircl{10}{0}\putbigcircl{20}{0}\end{picture}}$
and
${\setlength{\unitlength}{1.4pt}\begin{picture}(21,0)(5,-2)\thicklines\putbigsquar{10}{0}\putbigsquar{21}{0}\end{picture}}$,
but is not a union or join of coloured graphs (since we defined a
colouring as an equivalence relation).  Henceforth abbreviate
\emph{nicely coloured} to \emph{nice}.
\vspace{-.6ex}\begin{lemma}\label{fusion-nice}
A fusion of nice cographs is a nice cograph.
\end{lemma}\vspace{-2.5ex}
\begin{proof}
Let $C$ be the fusion of nice cographs $C_1$ and $C_2$
obtained by joining portions $C'_i$ of $C_i$.  Suppose $U$ is a union
of two-vertex colour classes in $C$ which induces a matching.  Let
$U_i=U\cap V(C_i)$ and $U'_i=U\cap V(C'_i)$.  By definition of fusion, the
only edges in $C$ between $U_1$ and $U_2$ are between $U_1'$ and
$U_2'$, and there are edges between all vertices of $U_1'$ and all
vertices of $U_2'$; thus ($\star$) \emph{there is at most one edge
between $U_1$ and $U_2$}, or else two edges of $C$ on $U$ would
intersect.  Since $U$ is a union of two-vertex colour classes, each
either in $U_1$ or $U_2$, each $U_i$ contains an even number of
vertices. Therefore, since $U$ induces a matching, ($\dagger$)
\emph{there must be an even number of edges between $U_1$ and $U_2$.}
Together ($\star$) and ($\dagger$) imply there is no edge between
$U_1$ and $U_2$, hence, for whichever $U_i$ is non-empty (perhaps
both), $U_i$ is a union of 
two-vertex
colour classes inducing a matching in
$C_i$, contradicting $C_i$ being nice.
\end{proof}
\vspace{-2.5ex}\begin{lemma}\label{fusion}
Every nice cograph with more than one colour class is a fusion of nice cographs.
\end{lemma}\vspace{-2.5ex}%
\begin{proof}
Let $C$ be a nice cograph.  Since $C$ is a cograph, its underlying
(uncoloured) graph has the form $(C_1\wedge C_2)\vee (C_3\wedge
C_4)\vee\!\ldots\!\vee(C_{n-1}\wedge C_n)\vee H$ for cographs $C_i$
and $H$ with no edges.  Assume $n\!\neq\! 0$, otherwise the result is
trivial.  Let $G$ be the graph whose vertices are the $C_i$,
with $C_iC_j\in E(G)$ iff there is an edge or colour class
$\{v,w\}$ in $C$ with $v\!\in
\!V(C_i)$ and $w\!\in\! V(C_j)$ (\textit{cf.}\ the proof of Theorem~4 in \cite{Ret03}).
A \emph{perfect matching} is a set of pairwise disjoint edges whose
union contains all vertices.  Since $C$ is nice, $M=\{
C_1C_2,C_3C_4,\ldots,C_{n-1}C_n\}$ is the only perfect matching of
$G$.  For if $M'$ is another perfect matching, then
$M'\mkern-2mu\setminus\mkern-2mu M$ determines a set of two-vertex
colour classes in $C$ whose union induces a matching in $C$: for each
$C_iC_j\in M'\mkern-2mu\setminus\mkern-2mu M$ pick a colour class
$\{v,w\}$ with $v\in V(C_i)$ and $w\in V(C_j)$.  Since $G$
has a unique perfect matching, some $C_kC_{k+1}\in M$ is a bridge
\cite{Kot59,LP86},
\textit{i.e.},
$(\,V(G),\,E(G){}\setminus{} C_kC_{k+1}\,)=X\vee
Y$ with $C_k\!\in\! V(X)$ and $C_{k+1}\!\in\! V(Y)$.  Let $W$ be the
union of all colour classes of $C$ coincident with any $C_i$ in
$X$,
and let $W'=V(C)\setminus W$. Then $C[W]$ and $C[W']$ are nice (since
$W$ and $W'$ are unions of colour classes), and $C$ is the fusion of
$C[W]$ and $C[W']$ joining portions $C_k$ of $C[W]$ and
$C_{k+1}$ of $C[W']$.
\end{proof}
\begin{lemma}\label{tensor-true}
Let $P_1$ and $P_2$ be combinatorial propositions and $Q$ a
combinatorial proposition or the empty graph.  Then
$(P_1\mkern-2mu\wedge\mkern-2mu P_2)\mkern-2mu\vee\mkern-2mu Q$ is
true iff $P_1\mkern-2mu\vee\mkern-2mu Q$ and
$P_2\mkern-2mu\vee\mkern-2mu Q$ are true.
\end{lemma}
\begin{proof}
A clause of $(P_1\mkern-2mu\wedge\mkern-2mu
P_2)\mkern-2mu\vee\mkern-2mu Q$ is a clause of $P_1\mkern-2mu\vee\mkern-2mu Q$
or $P_2\mkern-2mu\vee\mkern-2mu Q$, and vice versa.
\end{proof}
\vspace*{-2ex}\begin{theorem}[Combinatorial Soundness]\label{thm-soundness}
If a combinatorial proposition has a combinatorial proof, it is true.
\end{theorem}
\begin{proof}
Let $h:C\to P$ be a combinatorial proof.  We show $P$ is true by
induction on the number of colour classes in $C$.  In the base case,
$V(C)$ is a colour class.  If $v\in V(C)$ then $h(v)$ is in no
edge of $P$ (for if $h(v)w\in E(P)$ then a skew lifting at $v$ is an
edge in $C$, a contradiction), hence is in every clause $K$ of $P$.
Since $V(C)$ is axiomatic, $K$ is true.

\emph{Induction step.}  
By Lemmas~\ref{imgtrue} and \ref{shallow}, assume $h$ is shallow and
surjective.  By Lemma~\ref{fusion}, $C$ is a fusion of nice cographs
$\mkern1mu C_1$ and $C_2$ obtained from $C_1\mkern-3.5mu\vee C_2$ by
joining portions $C'_i$ of $C_i$.  If 
$C\mkern-2mu=\mkern-2muC_1\mkern-1.5mu\vee C_2$ then
$h'\mkern-1mu:\mkern-.5muC_1\mkern-1mu\to\mkern-1mu P$ defined by
$h'(v)\mkern-2mu=\mkern-2mu h(v)$ is a combinatorial proof, and $P$ is
true by induction hypothesis.  Otherwise each $C'_i$ is non-empty.
Let $P_i=h(C'_i)$.  Since $C'_1\wedge C'_2$ is a component of $C$ and
$h$ is a shallow surjection, $P_1\wedge P_2$ is a component of $P$,
say $P=(P_1\wedge P_2)\vee Q$.  Define $h_i:C_i\to P_i\vee Q$ by
$h_i(v)=h(v)$, a combinatorial proof: $C_i$ is a nice cograph, the
axiomatic colour class property is inherited from $h$, and $h_i$ is a
skew fibration by Lemma~\ref{long-lemma} (applied after forgetting
colourings).  By induction hypothesis $P_i\vee Q$ is true, hence $P$
is true by Lemma~\ref{tensor-true}.
\end{proof}
\vspace{-2ex}\begin{theorem}[Combinatorial Completeness]\label{abs-comp}
Every true combinatorial proposition has a combinatorial proof.
\end{theorem}
\noindent\textit{Proof.}
Let $P$ be a true combinatorial proposition.  We construct a
combinatorial proof of $P$ by induction on the number of edges in $P$.
In the base case $V(P)$ is a true clause, so there exists $W\subseteq
V(P)$ comprising a $1$-labelled vertex or a pair of vertices labelled
with dual literals.  Inclusion $W\to P$ is a combinatorial proof
(viewing $W$ as a graph with no edge and a single colour class, and
forgetting its labels).

\emph{Induction step.} 
Since $P$ is a cograph with an edge, $P=(P_1\wedge P_2)\vee Q$ for
combinatorial propositions $P_i$ and $Q$ a combinatorial proposition
or the empty graph.  Assume $Q$ is empty or not true; otherwise by
induction there is a combinatorial proof $C\to Q$ composable with
inclusion $Q\to P$ for a combinatorial proof of $P$, and we are done.
By Lemma~\ref{tensor-true}, $P_i\vee Q$ is true, so by induction has a
combinatorial proof $h_i:C_i\to P_i\vee Q$.  Let $C$ be
the fusion of $C_1$ and $C_2$ obtained by joining the portions
$h_i^{-1}(P_i)$ of $C_i$.  By Lemma~\ref{fusion-nice}, $C$ is nice.
Define $h:C\to P$ by $h(v)=h_i(v)$ iff $v\in V(C_i)$.  Then $h$ is a graph
homomorphism: let $vw\in E(C)$ with $v\in V(C_i)$ and $w\in V(C_j)$;
if $i\!=\!j$ then $h(v)h(w)\in E(P)$ since $h_i$ is a homomorphism; if
$i\neq j$ then $vw$ arose from fusion, so $h(v)\in P_i$ and $h(w)\in
P_j$, hence $h(v)h(w)\in E(P)$ since $P_1\wedge P_2\subseteq P$ has
all edges between $P_1$ and $P_2$.

The axiomatic colour class property for $h$ is inherited from the
$h_i$, so it remains to show that $h$ is a skew fibration.  Let $v\in
V(C)$ and $h(v)w\in E(P)$.  By symmetry, assume $v\in V(C_1)$.  Assume
$h(v)\!\in\! V(P_1)$ and $w\!\in\!V(P_2)$, otherwise we immediately
obtain a skew lifting of $h(v)w$ since $h_1$ is a skew fibration.
There is a vertex $x$ in $h_2^{-1}(P_2)$: if $Q$ is empty, this is
immediate; otherwise $Q$ is not true and $\restr{{h_2}}{Q}:C_2\to Q$
would be a combinatorial proof, contradicting soundness. Since fusion
joined the $h_i^{-1}(P_i)$, we have $vx\mkern-2mu\in\mkern-2mu E(C)$.
If $h(x)w\not\in E(P_2)$ we are done; otherwise since $h_2$ is a skew
fibration and $h(x)w\!\in\!  E( P_2)$ there exists $xy\!\in\! E(C_2)$
with $h(y)w\!\not\in\!  E(P_2)$.  Since $vy\!\in\! E(C)$ (again by
fusion), we have the desired skew lifting of $h(v)w$ at $v$.  (See
figure below.  Note $h(y)=w$ is possible.)
\begin{center}\vspace{5.5ex}
\newcommand{\qline}[4]{\qbezier(#1,#2)(#3,#4)(#3,#4)}
\small\setlength{\unitlength}{.55pt}\begin{picture}(0,65)(25,15)
\put(0,-8){
\putmathlabel{45}{122}{x}
\put(45,116){\vector(0,-1){51}}
\putmathlabel{100}{109}{y}
\put(99,101){\vector(0,-1){51}}
\putmathlabel{-40}{110}{v}
\put(-40,103){\vector(0,-1){51}}
}
\putmathlabel{45}{45}{h(x)}
\putmathlabel{55}{12}{w}
\putmathlabel{100}{30}{h(y)}
\putmathlabel{-40}{30}{h(v)}
\thicklines
\put(0,-8){\qline{-32}{112}{35}{121}
\qline{52}{120}{90.5}{112}
\qline{-32}{109}{90}{109}
}
\qline{-23.5}{32}{25}{42}
\qline{-23}{29}{81}{29}
\qline{-23.5}{26}{45}{14}
\qline{47}{36}{53}{19}
\qline{62}{40}{82}{33}
\putmathlabel{-41}{-1}{\underbrace{\hspace{4.3ex}}}
\putmathlabel{-40}{-16}{P_1}
\putmathlabel{71}{-1}{\underbrace{\hspace{11.4ex}}}
\putmathlabel{72}{-16}{P_2}
\end{picture}\vspace{.5ex}\end{center}
\hfill$\square$

\end{multicols}
\begin{twocolumn}
\newcommand{\authorfont}{}\end{twocolumn}
\begin{thebibliography}{\footnotesize BLS99\hspace{.2ex}}\vspace{-.4ex}\footnotesize\newcommand{\bibscrunch}{\vspace*{-.9ex}}
\bibitem[And81]{And81}
{\authorfont{}Andrews,~P.~B.}
\newblock{\em Theorem Proving via General Matings.}
\newblock{\sl J.\ ACM}\ 28 1981 193--214.\bibscrunch%
\bibitem[Bib74]{Bib74}
{\authorfont{}Bibel,~W.}
\newblock{\em An approach to a systematic theorem proving procedure in first-order logic.}
\newblock{\sl Computing}\ 12 1974 43--55.\bibscrunch%
\bibitem[BV02]{BV02}
{\authorfont{}Boldi,~P.~\& S.~Vigna.}
\newblock{\em Fibrations of graphs.}
\newblock{\sl Discrete Math.}\ 243 2002 21--66.\bibscrunch%
\bibitem[BLS99]{BLS99}
{\authorfont{}Brandst\"adt,~A., V.~B.~Le \& J.~P.~Spinrad.}
\newblock{\em Graph Classes: A Survey.}
\newblock{\sl SIAM monographs on Discr.\ Math.\ \& Applic.} 1999.\bibscrunch\bibitem[CR79]{CR79}
{\authorfont{}Cook,~S.~A.~\& R.~A.~Reckhow.}
\newblock{\em The relative efficiency of propositional proof systems.}
\newblock{\sl J.~Symb.~Logic}\/ 44(1) 1979 36--50.\bibscrunch\bibitem[CLS81]{CLS81}
{\authorfont{}Corneil,~D.~G., H.~Lerchs~\& L.~Stewart-Burlingham.}
\newblock{\em Complement reducible graphs.}
\newblock{\sl \mbox{Discr.\ Appl.}\ Math.}\ 3 1981 163--174.%
\bibscrunch\bibitem[Dav71]{Dav71}
{\authorfont{}Davydov,~G.~V.}
\newblock{\em The synthesis of the resolution method and the inverse method.}
\newblock{\sl Zapiski Nauchnykh Seminarov Lomi} 20 1971 24--35.
Translation in J.\ Sov.\ Math.\ 1(1) 1973 12--18.%
\bibscrunch\bibitem[Fr1879]{Fr1879}
{\authorfont{}Frege,~G.}
\newblock{\em Begriffsschrift, eine der arithmetischen nachgebildete Formelsprache des reinen Denkens.}
1879.  Translation
\cite[5--82]{vHe67}.\bibscrunch\bibitem[Gen35]{Gen35}
{\authorfont{}Gentzen,~G.}
\newblock{\em Untersuchungen \"uber das logische Schlie\ss en I, II.}
\newblock{\sl Mathematische Zeitschrift} 39 1935 176--210 405--431.
Translation \cite[68--131]{Gen69}.\bibscrunch\bibitem[Gen69]{Gen69}
{\authorfont{}Gentzen, G.}
{\em The Collected Papers of Gerhard Gentzen.}
Ed.\ M.\ E.\ Szabo.
North-Holland 1969.%
\newcommand{\com}{Comp.}\newcommand{\tcs}{Theor.\ \com\ Sci.}\bibscrunch\bibitem[Gir87]{Gir87}
{\authorfont{}Girard, J.-Y.}
\newblock{\em Linear logic}.
\newblock{\sl \tcs} 50 1987 1--102.\bibscrunch\bibitem[Gir91]{Gir91}
{\authorfont{}Girard, J.-Y.}
\newblock{\em A new constructive logic: Classical logic}.
\newblock{\sl Math.\ Struc.\ \com\ Sci.}\ 1(3) 1991 255-296.\bibscrunch\bibitem[Gra66]{Gra66}
{\authorfont{}Gray,~J.~W.}
{\em Fibred and cofibred categories.}
In {\sl Proc.\ Conf.\ on Categorical Algebra, La Jolla 1965}\/ 21--83. Springer 1966.\bibscrunch\bibitem[Gro59]{Gro59}
{\authorfont{}Grothendieck,~A.}
{\em Technique de descente et th\'eor\`emes d'exi\-stence en g\'eom\'etrie alg\'ebrique. I:
G\'en\'eralit\'es. Descente par mor\-phismes fid\`element plats.}\hspace*{-.1ex}
{\sl S\'eminaire Bourbaki}\hspace*{.1ex} 190 1959--60.\bibscrunch\bibitem[vHe67]{vHe67}
{\authorfont{}van Heijenoort, J.} (ed.).\ {\em From Frege to G\"odel. A Source Book
in Mathematical Logic 1879--1931.}  Harvard Univ.\ Press 1967.\bibscrunch\bibitem[Hil28]{Hil28}
{\authorfont{}Hilbert, D.}
\newblock{\em Die Grundlagen der Mathematik},
\newblock{\sl Abhandlungen aus dem Mathematischen Seminar der Hamburgischen Universit\"at} 6 1928 65-85.
Translation
\cite[465--479]{vHe67}.%
\bibscrunch\bibitem[Joh87]{Joh87}
{\authorfont{}Johnstone, P.~T.}
\newblock{\em Notes on Logic and Set Theory.}
Cambridge University Press 1987.\bibscrunch\bibitem[KM71]{KM71}
{\authorfont{}Kelly, G.~M.\ \& S.~Mac{ }Lane.}
\newblock{\em Coherence in closed categories.}
\newblock{\sl J.\ Pure Appl.\ Algebra} 1 1971 97--140; erratum, {\em ibid.}\ no.\ 2, 219.\bibscrunch\bibitem[KW84]{KW84}
{\authorfont{}Ketonen, J.\ \& R.~Weyhrauch.}
\newblock{\em A Decidable Fragment of Predicate Calculus.}
\newblock{\sl \tcs} 32 1984 297--307.\bibscrunch\bibitem[Kot59]{Kot59}
{\authorfont{}Kotzig, A.}
{\em On the theory of finite graphs with a linear factor.}
{\sl Mat.-Fyz.\ Casopis Slovensk. Akad. Vied}
9 \& 10 1959.%
\bibscrunch\bibitem[LP86]{LP86}
{\authorfont{}Lov\'asz, L.\ \& M. D. Plummer.}
{\em Matching Theory.} North-Holland 1986.%
\bibscrunch\bibitem[deM68]{deM68}
{\authorfont{}de~Morgan, A.}
\newblock{\em Review of J.\hspace*{-.15ex}~M.~\hspace*{-.2ex}Wilson's `Elementary Geometry'}.
\newblock{\sl The Athenaeum} 2 1868 71--73.\bibscrunch\bibitem[Pei58]{Pei58}
{\authorfont{}Peirce,~C.~S.}%
\newblock{\em Collected Papers.}
\newblock Ed.\ Hartshorne, Weiss \& Burks.  Harvard University Press 1931--5 \& 1958.%
\bibscrunch\bibitem[Ret03]{Ret03}
{\authorfont{}Retor\'e,~C.}
\newblock{\em Handsome proof-nets: perfect matchings and co\-graphs.}
\newblock{\sl \tcs}\ 294 2003 473--488.%
\bibscrunch\bibitem[Whi78]{Whi78}
{\authorfont{}Whitehead,~G.~\hspace{-.2ex}W.}
{\em Elements of Homotopy Theory.}
Springer 1978.\end{thebibliography}
\end{document}